\documentclass[12pt]{article}
\usepackage{amssymb}

\newcommand{\be}{\begin{equation}}
\newcommand{\ee}{\end{equation}}
\newcommand{\bea}{\begin{eqnarray}}
\newcommand{\eea}{\end{eqnarray}}
\newcommand{\barray}{\begin{array}}
\newcommand{\earray}{\end{array}}

\newcommand{\bitem}{\begin{itemize}}
\newcommand{\eitem}{\end{itemize}}
\newtheorem{teo}{Theorem}[section]
\newcommand{\bt}{\begin{teo}}
\newcommand{\et}{\end{teo}}
\newtheorem{Def}{Definition}[section]
\newcommand{\bd}{\begin{Def}}
\newcommand{\ed}{\end{Def}}
\newtheorem{lem}{Lemma}[section]
\newcommand{\bl}{\begin{lem}}
\newcommand{\el}{\end{lem}}
\newtheorem{prop}{Proposition}[section]
\newcommand{\bp}{\begin{prop}}
\newcommand{\ep}{\end{prop}}
\newtheorem{cor}{Corollary}[section]
\newcommand{\bc}{\begin{cor}}
\newcommand{\ec}{\end{cor}}
\newtheorem{ex}{Example}[section]
\newcommand{\bex}{\begin{ex}}
\newcommand{\eex}{\end{ex}}
\newtheorem{rem}{Remark}[section]
\newcommand{\br}{\begin{rem}}
\newcommand{\er}{\end{rem}}

\begin{document}

\begin{center}
{\Large \textbf{Duality in a special class of \\ submanifolds and
Frobenius manifolds\footnotetext[1]{The work was completed with the
financial support of the Russian Foundation for Basic Research
(grant no. 08-01-00054) and the Programme for Support of Leading
Scientific Schools (grant no. NSh-1824.2008.1).}}}
\end{center}

\smallskip

\begin{center}
{\large {O.I. Mokhov}}
\end{center}

\smallskip

Consider totally nonisotropic $N$-dimensional submanifolds $M^N$ in
$(N + L)$-di\-men\-si\-o\-nal pseudo-Euclidean spaces $\mathbb{E}^{N
+ L}_k$ ($N$-dimensional submanifolds that are not tangent to
isotropic cones of the ambient $(N + L)$-dimensional
pseudo-Euclidean space at their points).

Let $(z^1, \ldots, z^{N + L})$ be pseudo-Euclidean coordinates in
$\mathbb{E}^{N + L}_k$; a submanifold $M^N$ is given by a smooth
vector-function $r = (z^1 (u^1, \ldots, u^N), \ldots, z^{N + L}
(u^1, \ldots, u^N))$, ${\rm rank\,} (\partial z^i /
\partial u^j) = N,$ $1 \leq i \leq N + L,$ $1 \leq j \leq N$;
$\partial r / \partial u^i = r_i (u)$, $1 \leq i
\leq N$, is a basis of the tangent space at an arbitrary point $u =
(u^1, \ldots, u^N)$ of the submanifold $M^N$; $n_1 (u), \ldots, n_L
(u)$ is an arbitrary basis of the normal space $\mathbb{N}_u$ that
depends smoothly on the point $u$; $g_{ij} (u) = (r_i, r_j)$, $1
\leq i, j \leq N$, is the first fundamental form of the submanifold
($(\cdot, \cdot)$ is the pseudo-Euclidean scalar product in
$\mathbb{E}^{N + L}_k$), $h_{\alpha \beta} (u) = (n_{\alpha},
n_{\beta})$, $1 \leq \alpha, \beta \leq L$ ($\det g_{ij} (u) \neq 0$
and $\det h_{\alpha \beta} (u) \neq 0$ for totally nonisotropic
submanifolds). The Gauss and Weingarten decompositions have the form
 \be
{\partial^2 r \over
\partial u^i \partial u^j} = a^k_{ij} (u) {\partial r \over \partial
u^k} + b^{\beta}_{ij} (u) n_{\beta} (u), \ \ \ \ \ \ \ {\partial
n_{\alpha} \over \partial u^j} = c^k_{\alpha j} (u) {\partial r
\over \partial u^k} + d^{\beta}_{\alpha j} (u) n_{\beta} (u),
\label{1} \ee respectively. The coefficients $a^k_{ij} (u),
b^{\beta}_{ij} (u), c^k_{\alpha j} (u), d^{\beta}_{\alpha j} (u)$,
the metric $g_{ij} (u)$ and the functions $h_{\alpha \beta} (u)$
satisfy a number of relations including the Gauss equations, the
Codazzi equations and the Ricci equations for any submanifold. Note
that in this paper we consider only the local theory of
submanifolds.

We shall single out a special class of $N$-dimensional submanifolds
in $2N$-di\-men\-si\-o\-nal pseudo-Euclidean spaces. Our purpose is
to single out the case, when the sets of the basis vectors of the
tangent and normal spaces, $r_i (u)$, $1 \leq i \leq N$, and
$n_{\alpha} (u)$, $1 \leq \alpha \leq L$, possess equal rights and
are dual one to another. It is obvious that the condition $L = N$ is
necessary for such duality. Moreover, another necessary condition is
the potentiality of the basis $n_{\alpha} (u)$: $n_{\alpha} (u) =
\partial n / \partial u^{\alpha}$, $1 \leq \alpha \leq N$, where $n
(u)$ is a certain vector-function on the submanifold.

{\bf Definition 1.} A basis $n_{\alpha} (u)$, $1 \leq \alpha \leq
N$, in the normal space $\mathbb{N}_u$ of a certain $N$-dimensional
submanifold in a $2N$-dimensional pseudo-Euclidean space is called
{\it potential} if there exists a vector-function (a {\it potential}
of the normals) $n (u)$ on the submanifold such that $n_{\alpha} (u)
= \partial n / \partial u^{\alpha}$, $1 \leq \alpha \leq N$.

If the vector-function $n (u)$ exists, then it generates a potential
basis in the normal space $\mathbb{N}_u$ in any coordinate system.

{\bf Definition 2.} An $N$-dimensional submanifold in a
$2N$-dimensional pseudo-Euclidean space is called a {\it submanifold
with potential normal basis} (or a {\it submanifold with potential
of normals}) if there exists a vector-function (a {\it potential of
normals}) $n (u)$ on the submanifold such that the vectors $n_i (u)
=
\partial n /
\partial u^i$, $1 \leq i \leq N$, form a basis in the normal space $\mathbb{N}_u$
at any point $u$ on the submanifold.

Definition 2 is invariant (it does not depend on a local coordinate
system).

The coefficients $a^k_{ij} (u)$ in the Gauss decomposition (1) are
always the coefficients of the Levi-Civita connection of the metric
$g_{ij} (u)$. If the submanifold is equipped with a potential normal
basis, then a similar assertion is also true for the coefficients
$d^{\beta}_{\alpha j} (u)$ in the Weingarten decomposition (1).
Consider the corresponding Gauss and Weingarten decompositions: \be
{\partial^2 r \over
\partial u^i \partial u^j} = a^k_{ij} (u) {\partial r \over
\partial u^k} + b^k_{ij} (u) {\partial n \over \partial u^k}, \ \ \
\ \ \ \ {\partial^2 n \over \partial u^i \partial u^j} = c^k_{i j}
(u) {\partial r \over \partial u^k} + d^k_{i j} (u) {\partial n
\over \partial u^k}. \label{2} \ee

{\bf Theorem 1.} {\it The functions $h_{ij} (u) = (\partial
n/\partial u^i,
\partial n/\partial u^j)$ define a covariant metric on the submanifold,
and the coefficients $d^k_{i j} (u)$ in the Weingarten decomposition
{\rm (2)} are the coefficients of the symmetric affine connection
compatible with the metric $h_{ij} (u)$, that is the coefficients of
the Levi-Civita connection of the metric $h_{ij} (u)$. The
coefficients $b^k_{ij} (u)$ and $c^k_{i j} (u)$ are tensors of type
{\rm (1, 2)} symmetric with respect to the subscripts $i$ and $j$ on
the submanifold $M^N$.}

{\bf Theorem 2 (duality principle).} {\it If a submanifold is
equipped with a potential of normals $n (u)$ and given by the
vector-functions $(r (u), n (u))$, then the vector-functions $(n
(u), r (u))$ also give a submanifold equipped with the potential of
normals $r (u)$ such that $\partial n/\partial u^i,$ $1 \leq i \leq
N,$ are tangent vectors and $\partial r/\partial u^i,$ $1 \leq i
\leq N,$ are basis normal vectors of the submanifold, and moreover,
in this case all the objects of the local theory of such
submanifolds are dual to each other, in particular, the Gauss
decomposition becomes the Weingarten decomposition and the
Weingarten decomposition becomes the Gauss decomposition, the Gauss
equations becomes the Ricci equations and the Ricci equations
becomes the Gauss equations, the Codazzi equations changes to
themselves {\rm (}they are self-dual{\rm )}, the tensor $b^k_{ij}
(u)$ becomes the tensor $c^k_{ij} (u)$ and the tensor $c^k_{ij} (u)$
becomes the tensor $b^k_{ij} (u)$, the metric $g_{ij} (u)$ becomes
the metric $h_{ij} (u)$ and vice versa.}

Submanifolds with potential of normals form an important and rich
class of submanifolds. Arbitrary one-dimensional submanifolds of
pseudo-Euclidean planes are a trivial example of such submanifolds.
The general theory of submanifolds with potential of normals, the
duality principle for them and important examples will be presented
in a separate paper. The submanifolds equipped with natural
Frobenius structures, which were constructed by the present author
in [1]--[5] and which realize arbitrary Frobenius manifolds (the
theory of Frobenius manifolds was constructed in [6]), are a
particular case of submanifolds with potential of normals. The
present author has proved in [1] and [2] that an arbitrary Frobenius
manifold can be realized as a certain flat submanifold with
potential of normals for which $g_{ij} (u) = c h_{ij} (u)$, $c =
{\rm const} \neq 0$, where $c$ is a deformation parameter preserving
the corresponding Frobenius structure. Let $h_{ij} (u) = \eta_{ij}$,
$\eta_{ij} =\eta_{ji}$, $\det \eta_{ij} \neq 0$, $\eta_{ij} = {\rm
const}$, $g_{ij} = c \eta_{ij}$, $c = {\rm const} \neq 0$. In this
case, for the submanifolds with potential of normals, the relations
$a^k_{ij} (u) = d^k_{i j} (u) = 0$ are satisfied, and also there
exists a function $\Phi (u)$ such that $b^k_{ij} (u) = \eta^{ks}
\partial^3 \Phi / \partial u^s \partial u^i \partial u^j,$
$c^k_{i j} (u) = - (1/c)\, b^k_{ij} (u)$, $\eta^{is} \eta_{sj} =
\delta^i_j$, and all the relations of the local theory of
submanifolds are satisfied if and only if the function $\Phi (u)$
satisfies the associativity equations of two-dimensional topological
quantum field theories (the Witten--Dijkgraaf--Verlinde--Verlinde
equations, see [6]) \be {\partial^3 \Phi \over
\partial u^i \partial u^j \partial u^s} \eta^{sp} {\partial^3 \Phi
\over \partial u^p \partial u^k \partial u^l} = {\partial^3 \Phi
\over \partial u^i \partial u^k \partial u^s} \eta^{sp} {\partial^3
\Phi \over \partial u^p \partial u^j
\partial u^l}. \ee

\smallskip

\begin{center}
{\bf Bibliography}
\end{center}

\smallskip

\noindent [1] O.I. Mokhov, Frobenius manifolds as a special class of
submanifolds in pseudo-Euclidean spaces, will be published in AMS
Translations Series 2, Volume 224, Advances in Mathematical
Sciences, eds. V.M.Buchstaber and I.M.Krichever, ``Geometry,
Topology, and Mathematical Physics: S.P.Novikov's Seminar:
2006--2007'', Amer. Math. Soc., Providence, RI 2008;
http://arXiv.org/abs/0710.5860 (2007).

\noindent [2] O.I. Mokhov, Theory of submanifolds, the associativity
equations in 2D topological quantum field theories, and Frobenius
manifolds, In: Proceedings of the Workshop ``Nonlinear Physics.
Theory and Experiment. IV.'' (Gallipoli, Lecce, Italy, 22 June -- 1
July, 2006), published in Teoret. i Matemat. Fizika {\bf 152}:2
(2007), 368--376; English transl., Theoret. Math. Phys. {\bf 152}:2
(2007), 1183--1190; Preprint MPI 06-152, Max-Planck-Institut f\"ur
Mathematik, Bonn, Germany, 2006; arXiv:math.DG/0610933 (2006).

\noindent [3] O.I. Mokhov, Submanifolds in pseudo-Euclidean spaces
and Dubovin--Frobenius structures, In: Proceedings of the 10th
International Confrence ``Differential Geometry and its
Applications'' in honour of the 300th anniversary of the birth of
Leonhard Euler, August 27--31, 2007, Olomouc, Czech Republic, World
Scientific, Singapore, 2008, pp. 505--516.

\noindent [4] O.I. Mokhov, Non-local Hamiltonian operators of
hydrodynamic type with flat metrics, and the associativity
equations, Uspekhi Mat. Nauk {\bf 59}:1 (2004), 187--188; English
transl., Russian Math. Surveys {\bf 59}:1 (2004), 191--192.

\noindent [5] O.I. Mokhov, Nonlocal Hamiltonian operators of
hydrodynamic type with flat metrics, integrable hierarchies, and the
associativity equations, Funktsional. Anal. i Prilozhen. {\bf 40}:1
(2006), 14--29; English transl., Funct. Anal. Appl. {\bf 40}:1
(2006), 11--23; arXiv:hep-th/0406292 (2004).

\noindent [6] B. Dubrovin, Geometry of 2D topological field
theories, In: Lecture Notes in Math., vol. 1620, Springer-Verlag,
Berlin, 1996, pp. 120--348; arXiv:hep-th/9407018 (1994).

\begin{flushleft}
{\bf O.I. Mokhov}\\
Centre for Non-Linear Studies,\\
Landau Institute for Theoretical Physics,\\
Russian Academy of Sciences;\\
Department of Geometry and Topology,\\
Faculty of Mechanics and Mathematics,\\
Lomonosov Moscow State University\\
{\it E-mail\,}: mokhov@mi.ras.ru; mokhov@landau.ac.ru; mokhov@bk.ru\\
\end{flushleft}

\end{document}